\documentclass[11pt]{article}
\usepackage{latexsym, amsmath, amssymb}
\usepackage{amsmath}
\usepackage{amsfonts}
\usepackage{amssymb}
\usepackage{amscd}
\usepackage{amsthm}
\usepackage{fancyhdr}
\newtheorem{theorem}{Theorem}[section]

\newtheorem{remark}[theorem]{Remark}
\newtheorem{lemma}[theorem]{Lemma}
\newtheorem{definition}[theorem]{Definition}
\newtheorem{corollary}[theorem]{Corollary}

\newcommand{\rk}{\mathrm{rank}}

\def\Z{\mathbb Z}

\def\C{\mathbb C}

\def\Q{\mathbb Q}

\newcommand\sE{{\mathcal E}}

\newcommand\sI{{\mathcal I}}

\newcommand\sO{{\mathcal O}}

\newcommand\sS{{\mathcal S}}

\newcommand\sZ{{\mathcal Z}}

  \def\peen{\hbox{$ {\mathbf  P}^n$}}

  \def \tab#1{\kern #1 truein}

  \def\E{\hbox{${\cal E}$}}
  \def\F{\hbox{${\cal F}$}}
  \def\G{\hbox{${\cal G}$}}
  \def\A{\hbox{${{\cal A}}$}}
  \def\B{\hbox{${\cal B}$}}
  \def\C{\hbox{${\cal C}$}}

  \def\Q{\hbox{${\cal Q}$}}

\begin{document}
\title{Monads and Vector Bundles on Quadrics} 
\date{}
\author{Francesco Malaspina \\
 Dipartimento di Matematica Universit\`a di Torino\\
via Carlo Alberto 10, 10123 Torino, Italy\\ 
{\small\it e-mail: francesco.malaspina@unito.it}}

   \maketitle%\maketitle\footnotetext{Written with the support of the University Ministry funds.}
   \def\thefootnote{}
  \baselineskip=+.5cm\begin{abstract}We improve Ottaviani's splitting criterion for vector bundles on a quadric hypersurface and obtain the equivalent of the result by Rao, Mohan Kumar and Peterson. Then we give the classification of rank $2$ bundles without "inner" cohomology on $\Q_n$ ($n>3$). It surprisingly  exactly agrees with the classification by Ancona, Peternell and Wisniewski of rank $2$ Fano bundles.
      \end{abstract}

   \section*{Introduction}

A monad on $\mathbb{P}^n$ or, more generally, on a projective variety $X$, is a complex of three vector bundles
$$0 \rightarrow \A \xrightarrow{\alpha} \B \xrightarrow{\beta} \C \rightarrow 0$$

such that $\alpha$ is injective and $\beta$ is surjective.
Monads have been studied by Horrocks, who proved (see \cite{Ho} or \cite{BH}) that every vector bundle on $\mathbb{P}^n$  is the homology of a suitable minimal monad.
Throughout the paper we often use the Horrocks correspondence between a bundle $\E$ on $\mathbb P^n$ ($n\geq 3$) and the corresponding minimal monad
  $$0 \rightarrow \A \xrightarrow{\alpha} \B
\xrightarrow{\beta} \C \rightarrow 0,$$ where $\A$ and $\C$ are sums of line
bundles and $\B$ satisfies:
\begin{enumerate}
\item $H^1_*(\B)=H^{n-1}_*(\B)=0$ 
\item $H^i_*(\B)=H^i_*(\E)$ \ $\forall i, 1<i<n-1
$.
\end{enumerate}
This correspondence holds also on $X$ ($\dim X\geq 3$). Indeed the proof of the result in (\cite{BH} proposition $3$) can be easily extended to $X$ (see \cite{Ml}) theorem 2.1.6.). \\
Rao, Mohan Kumar and Peterson  have successfully used this tool to investigate the intermediate cohomology modules of a vector bundle on $\mathbb{P}^n$ and give cohomological splitting conditions (see \cite{KPR1}).\\
The first aim of the present paper is to extend to smooth quadric hypersurfaces the above result by Rao, Mohan Kumar and Peterson. In $\Q_n$ the Horrocks criterion does not work, but there is a theorem that classifies all the ACM bundles (see \cite{Kn}) as  direct sums of line bundles and spinor bundles (up to a twist - for generalities about spinor bundles see \cite{Ot2}).\\

In the first section we prove some necessary conditions that  a minimal monad associated to a bundle $\sE$ must satisfy.\\
%\newpage
The second aim of this paper is the improvement of Ottaviani's splitting criterion  (see \cite{Ot1} and \cite{Ot3}): we obtain the equivalent of the result by Rao, Mohan Kumar and Peterson  on a quadric hypersurface. 
In the last section we focus our interest on rank two vector bundles  on $\Q_4$ and prove the following theorem, which is  our main result: 

For an  indecomposable rank $2$ bundle $\E$ on $\Q_4$ with $H_*^1(\E)\not=0$ and $H_*^2(\E)=0$, the only possible minimal monad, such that both $\A$ and $\C$ do not vanish,  is (up to a twist)
\begin{equation} 0 \to \sO \rightarrow
\sS'(1)\oplus\sS''(1) \rightarrow \sO(1) \to 0,
\end{equation}
and such a monad exists.

This means that the two spinor bundles  and the bundle corresponding to this monad  are the only rank $2$ bundles without "inner" cohomology (i.e. $H^2_*(\E)= ...=H^{n-2}_*(\E)=0$).
 By using monads again we can also understand the behavior of rank two bundles on $Q_5$  and also on $\Q_n$, $n >5$. More precisely we can prove that:
\begin{enumerate}
\item For an  indecomposable rank $2$ bundle $\E$ on $\Q_5$ with
$H_*^2(\E)=0$ and $H_*^3(\E)=0$, the only possible minimal monad, such that both $\A$ and $\C$ do not vanish,  is (up to a twist)
$$ 0 \to \sO \rightarrow
\sS_5(1) \rightarrow \sO(1) \to 0,$$
and such a monad exists.
 \item For $n>5$, there is no indecomposable bundle of rank $2$ on $\Q_n$ with\\
$H^2_*(\E)= ...=H^{n-2}_*(\E)=0$.
\end{enumerate}
It is surprising that this classification of rank $2$ bundle on $\mathbb P^n$ and $\Q_n$ ($n>3$) exactly agrees with the classification by Ancona, Peternell and Wisniewski of rank $2$ Fano bundles (see \cite{APW}).\\
We can say that if $\E$ is a rank $2$ bundle on $\mathbb P^n$ and $\Q_n$ ($n>3$), then $$\textrm{$\E$ is a Fano bundle $\Leftrightarrow \E$ is without inner cohomology.}$$
I would like to thank A. Prabhakar Rao for having introduced me into the topic and Giorgio Ottaviani for his useful comments and suggestions.\\

\section{Monads for Bundles without inner cohomology}
In this section $X$ denotes a nonsingular subcanonical, irreducible ACM projective variety.\\
If $M$ is a finitely generated module over the homogeneous coordinate ring of $X$, we denote by  $\beta_i(M)$ the total Betti numbers of $M$.\\
We say that a bundle is indecomposable if it does not split as a direct sum of line bundles.\\ 
\begin{definition}We will call bundle without inner cohomology a bundle $\E$ on $X$ with 
$$H^2_*(\E)= \dots =H^{n-2}_*(\E)=0,$$
where $n=dim X$.
\end{definition}
 In $\mathbb P^n$ Kumar Peterson and Rao showed that,
if $n$ is even and  $\textrm{$\rk$}(\mathcal E)< n$ (or if $n$ is odd and 
$\textrm{$\rk$}(\mathcal E)< n-1$),
and $$ 0 \to \A \xrightarrow{\alpha} \B
\xrightarrow{\beta} \C \to 0$$ is a minimal monad for $\E$ such that $\A, \B$ and $\C$ are not
zero, then $\B$ cannot split.\\
This means that $\mathcal E$ splits if and only if it is without inner cohomology.

On $X$ we are able to prove  the first part of the theorem about monads:\\
   \begin{theorem}\label{t3}Let $\mathcal E$ be a vector bundle on  $X$
     of dimension $n$, with $n>3$.
\begin{enumerate}
\item \label{firsts}  If $n$ is even and if $$\textrm{$\rk$}(\mathcal E)< n$$
then no minimal monad for $\E$ exists such that $\A$ or $\C$
are not zero and $\B$ is split. \item \label{seconds} If $n$ is
odd and if $$\textrm{$\rk$}(\mathcal E)< n-1$$ then no
minimal monad for $\E$ exists such that $\A$ or $\C$ are not zero and $\B$ is
split.
\end{enumerate}\end{theorem}
First of all we prove a simple and useful lemma:\\

\begin{lemma} Let $\E$ be a bundle on $X$  with $H_*^2(\E)=H^{n-2}_*(\E)=0$ where $n=\dim X>3$ and let $H$ be  a hyperplane such that $X'=X\cap H$ is a subcanonical, irreducible, ACM, nonsingular projective variety. (use Bertini's theorem for irreducibility).\\
If 
$$ 0 \to \A \rightarrow
\B \rightarrow \C \to 0$$  is a minimal monad for $\E$, then a minimal monad for the restriction
$\E_{|X'}$  is just the restriction
$$ 0 \to \A_{|X'} \rightarrow
\B_{|X'} \rightarrow \C_{|X'} \to 0.$$
\begin{proof}From the sequence
$$0 \to \E(-1) \rightarrow \E \rightarrow \E_{|X'} \to 0,$$
and the corresponding sequence in cohomology
$$H^1_*(\E(-1)) \rightarrow H^1_*(\E) \xrightarrow{\gamma} H^1_*(\E_{|X'})
\to H^2_*(\E(-1))=0$$ we see that the map $\gamma$ is surjective.  Then the module $H^1_*(\E_{|X'})$ has
the same generators of $H^1_*(\E)$ of the same degrees restricted to $X'$ and this
means that, if $$ 0 \to \A' \xrightarrow{\alpha} \B'
\xrightarrow{\beta} \C' \to 0$$ is a minimal monad for $\E_{|X'}$, then $$\C'\simeq\C_{X'}.$$

In the same way, by using the fact that $H^{n-2} _*(\E)=0$, we see that $$\A'\simeq \A_{|X'}.$$  Then, by construction we see that
 also
$$\B'\simeq\B_{|X'}.$$
\end{proof}\label{l1}
\end{lemma}

\begin{proof} (of theorem \ref{t3})

Let us suppose that we know the result of the theorem for $n$ even. Let \E\ be
a bundle on $X$ with $$\textrm{$\rk$}(\mathcal E) < n-1,$$ $n> 3$, $n$ odd.
Let us also suppose that we have a minimal monad $$ 0 \to \A \xrightarrow{\alpha} \B
\xrightarrow{\beta} \C \to 0,$$ where $\A$ and $\C$ not zero and
$\B$ splits.\\
 Let $H$ be any hyperplane such that $X'=X\cap H$ is a subcanonical, irreducible, ACM projective variety. By (\ref{l1}) we have that $$ 0 \to
\A_{|X'} \rightarrow \B_{|X'} \rightarrow \C_{|X'} \to 0$$ is the
minimal monad for $\E_{|X'}$, where $\E_{|X'}$ is a bundle of rank $<n-1$, and 
$n-1=\dim X$ is even.\\
 Now if $\B$ splits also $\B_{|X'}$ has
to split and this is against our assumption of the result of the
theorem for $n$ even.\\
 Thus, establishing the result of the
theorem for the case of $n$ even will also establish the result
for $n$ odd.

 Now if one of $\A$ or $\C$ is zero, and $\B$ splits, then either $\E$ or its dual is a first
syzygy module. In this case, \E\ must have rank at least $n$ by
the following argument.\\
 Assume that \C\ is zero and let
$r$ be the rank of \E. From the short exact sequence $$0 \to \A
\to \B \to \E \to 0$$ we get the exact sequence $$ 0 \rightarrow S^r\A
\xrightarrow{\gamma_{r} } S^{r-1}\A\otimes \wedge^1\B \xrightarrow{\gamma_{r-1}}\dots \xrightarrow{\gamma_{2}}
 S^{1}\A\otimes \wedge^{r-1}\B \xrightarrow{\gamma_{1}}\wedge^r\B \xrightarrow{\gamma_{0}}
  \wedge^r\E \rightarrow 0.$$
  If we put $$\Gamma_i=ker\gamma_i,$$ we see that $$H^r_*(S^r\A)=H^r_*(\Gamma_{r-1})=\dots =H^i_*(\Gamma_{i-1})=\dots =H^1_*(\Gamma_0)$$ and $\forall 0< j<r$,
  $$H^{r-j}_*(S^r\A)=H^{r-j}_*(\Gamma_{r-1})=\dots =H^{i-j}_*(\Gamma_{i-1})=\dots =H^1_*(\Gamma_j).$$
When $r < n$, $$H^r_*(S^r\A)=0,$$ so $$H^1_*(\Gamma_0)=0$$ and  $$H^0_*(\wedge^r\B) \rightarrow
H^0_*(  \wedge^r\E)$$ is a surjective map between free modules.\\
This means that the map $\gamma_0$ splits and the bundle $\Gamma_0$ is a direct sum of line bundles.\\
Now, since also $$H^{r-1}_*(S^r\A)=\dots =H^1_*(S^r\A)=0,$$
we have $$H^1_*(\Gamma_1)=\dots =H^1_*(\Gamma_{r-1})=0.$$
We consider, then, the short exact sequence $$0\rightarrow \Gamma_1 \rightarrow S^{1}\A\otimes \wedge^{r-1}\B \rightarrow \Gamma_0\rightarrow 0.$$ Since $\Gamma_0$ is free and $$H^1_*(\Gamma_1)=0$$ we have that also this sequence splits and, hence, the map $\gamma_1$ splits and the bundle $\Gamma_1$ is a sum of line bundles.\\ 
 By iterating this argument we can conclude that the long
exact sequence is split at each place.\\ In particular, the map
$$S^r\A \to S^{r-1}\A\otimes \wedge^1\B,$$  which is obtained from
$\alpha$ as $$a_1a_2\dots a_r \to \sum (\pm a_1a_2..\hat a_i..a_r
\otimes \alpha(a_i)),$$ is split.\\ This goes against the minimality of the monad.\\
Suppose now that $\A$ and $\C$ are both not zero and $n$ is even
with $n=2k$.\\
 Let \E\ be a bundle on $X$ with $$\textrm{$\rk$}(\mathcal E)
\leq n-1.$$ By adding line bundles to \E\ (if necessary), we may
suppose that $$\textrm{$\rk$}(\mathcal E) = n-1.$$ 
Now we can follow the proof in (\cite{KPR1} pages 7-8] and see that such a monad \[0 \to \A \xrightarrow{\alpha} \B
\xrightarrow{\beta} \C \to 0\]
cannot exist.

 \end{proof} 
 \begin{remark}\label{r1}The Kumar-Peterson-Rao  theorem tells us that on $\mathbb P^n$ there is no indecomposable bundle without inner cohomology with small rank.

In a more general space $X$ we cannot say that because the Horrocks theorem fails.
But is still true the following:\\
In a minimal monad 
\[0 \to \A \xrightarrow{\alpha} \B
\xrightarrow{\beta} \C \to 0,\] for a bundle without inner cohomology on $X$, the bundle $\B$ must be ACM and indecomposable.
\end{remark}

 Now we prove a theorem about minimal monads for rank $2$ bundles: 

\begin{theorem}\label{t2} Let $X$ be 
 of dimension $n>3$, and $\E$ a rank $2$ bundle with $H^2_*(\E)=0$. Then any minimal monad $$ 0 \to \A \xrightarrow{\alpha} \B
\xrightarrow{\beta} \C \to 0$$ for $\E$, such that $\A, \B$ and $\C$ are not
zero, must satisfy the
following conditions:
\begin{enumerate}

\item $H^1_*(\wedge^2\B)\not=0$ and $ \beta_0(H^1_*(\wedge^2\B))\geq
\beta_0(H^0_*(S_2\C)).$
 \item $H^2_*(\wedge^2\B)=0$
\end{enumerate}
\begin{proof}First of all, since $X$ is ACM, the sheaf $\sO_X$ does not have intermediate
cohomology. The same is true for $\A$ and $\C$ that are free $\sO_X$-modules.\\
Let us now  assume the existence of a minimal monad with $H^1_*(\wedge^2\B)=0$
$$ 0 \to \A \xrightarrow{\alpha} \B
\xrightarrow{\beta} \C \to 0.$$ Then, if $\G=\ker\beta$,
from the sequence
$$ 0 \to S_2\A  \to (\A \otimes
\G)\to \wedge^2 \G
   \to \wedge^2\mathcal E \to 0, $$
   we have
$$ H^2_*(\wedge^2 \G) = H^2_*(\A \otimes \G) = 0,$$ since $H^2_*(\B)=H^2_*(\G)=H^2_*(\E)=0$.\\
Moreover, from the sequence
$$ 0 \to \wedge^2\G \to \wedge^2\B \to \B\otimes \C
   \to S_2\C \to 0, $$
 passing to the exact
sequence of maps on cohomology groups, since $H^1_*(\wedge^2\B)=H^2_*(\wedge^2 \G)=0$, we get
$$    H^0_*(\B\otimes \C)
   \to H^0_*(S_2\C) \to 0 .$$
Now, if we call $S_X$ the coordinate ring, we can say that $H^0_*(S_2\C)$
   is a free $S_X$-module, hence projective; then there exists a map $$H^0_*(\B\otimes \C)
   \leftarrow H^0_*(S_2\C)$$ and this means that $$\B\otimes \C
   \to S_2\C \to 0 $$ splits.\\
   But this map is obtained from $\beta$ as $b\otimes c\mapsto
   \beta(b)c$, so if it splits also $\beta$ has  to split and this violates the
   minimality of the monad.
We can say something stronger.\\
From the sequence
$$0 \to \wedge^2\G \to \wedge^2\B \to \B\otimes \C\xrightarrow{\gamma} S_2\C \to 0, $$
    since $H^2_*(\wedge^2 \G)  = 0$, we have a surjective map
    $$    H^1_*(\wedge^2\B)
   \to H^1_*(\Gamma) \to 0 $$ where $\Gamma=\ker\gamma$, and then
   $$     \beta_0(H^1_*(\wedge^2\B))\geq \beta_0(H^1_*(\Gamma)).$$
   On the other hand we have the sequence

    $$    H^0_*(\B\otimes \C)
   \xrightarrow{\gamma} H^0_*(S_2\C)\rightarrow  H^1_*(\Gamma) \to 0; $$ so, if
    $$     \beta_0(H^1_*(\wedge^2\B))< \beta_0(H^0_*(S_2\C)),$$ also $$\beta_0(H^1_*(\Gamma))< \beta_0(H^0_*(S_2\C)),$$ and some of the generators
    of $H^0_*(S_2\C)$ must be in the image of $\gamma$.\\
    But $\gamma$ is obtained from $\beta$ as $b\otimes c\mapsto
   \beta(b)c$, so  also  some generators of $C$ must be in the image of $\beta$ and this contradicts the
   minimality of the monad.\\
   We conclude that not just $H^1_*(\wedge^2\B)$ has to be non zero
   but also $$     \beta_0(H^1_*(\wedge^2\B))\geq \beta_0(H^0_*(S_2\C)).$$
The second condition
comes from the sequence $$ 0 \to \wedge^2\G \to \wedge^2\B \to
\B\otimes \C
   \to S_2\C \to 0, $$ since $H^2_*(\wedge^2\G)=H^2_*(
\B\otimes \C)=0$.
   \end{proof}
   \end{theorem}
   %\begin{remark}If $X$ of dimension $n$ with $n>3$, there is not minimal monads  for a rank $2$ bundle $$ 0 \to \A \xrightarrow{\alpha} \B
%\xrightarrow{\beta} \C \to 0,$$ such that $\A$ and $\C$ are not zero
%and $\B$ is split.
%\end{remark}

%%%%%%%%%%%%%%%%%
%%%%%%%%%%%%%%%%%

\section{Splitting Criteria on $\Q_n$}
In this section we apply our results to a smooth quadric hypersurface $\Q_n$ in $\mathbb P^{n+1}$.\\
Let us notice that $\Q_n$ is a nonsingular, ACM, irreducible projective variety and, if $n>3$, we  also have $$\textrm{Pic}(\Q_n)=\mathbb Z,$$ so it satisfies all the conditions of $X$.\\
First of all we need a useful remark about spinor bundles:\\

\begin{remark}\label{r2}
By applying (\cite{Ot2} Lemma 2.7. and Theorem 2.8)  we have that if $n=2m+1$, 
$$h^0(\Q_n,\sS(1)\otimes\sS)=1.$$ So from the sequence
$$ 0 \to \sS\otimes \sS
\rightarrow \sO_{\Q_n}^{2^{m+1}}\otimes \sS \rightarrow \sS(1)\otimes \sS \to 0,$$ and the sequence in cohomology $$ 0 =H^0( \sO_{\Q_n}^{2^{m+1}}\otimes \sS) \rightarrow H^0(\sS(1)\otimes \sS)\rightarrow H^1(\sS\otimes\sS)\to 0$$ we see that $$h^0(\Q_n,\sS(1)\otimes\sS)=h^1(\Q_n,\sS\otimes\sS)=1.$$
Moreover, if $n=4m$, we have $$h^0(\Q_n,\sS'(1)\otimes\sS')=h^0(\Q_n,\sS''(1)\otimes\sS'')=1$$ and $$h^0(\Q_n,\sS'(1)\otimes\sS'')=h^0(\Q_n,\sS''(1)\otimes\sS')=0.$$ Then 
$$h^1(\Q_n,\sS''\otimes\sS')=1$$ and $$h^1(\Q_n,\sS''\otimes\sS'')=h^1(\Q_n,\sS'\otimes\sS')=0$$

 while, if $n=4m+2$, $$h^0(\Q_n,\sS'(1)\otimes\sS')=h^0(\Q_n,\sS''(1)\otimes\sS'')=0$$ and $$h^0(\Q_n,\sS'(1)\otimes\sS'')=h^0(\Q_n,\sS''(1)\otimes\sS')=1.$$ Then 
 $$h^1(\Q_n,\sS''\otimes\sS')=0$$ and $$h^1(\Q_n,\sS''\otimes\sS'')=h^1(\Q_n,\sS'\otimes\sS')=1.$$
\end{remark}
Our starting point is the splitting criterion of Ottaviani (see \cite{Ot1} or \cite{Ot3}).
By using monads we can improve this criterion in the case of bundle with a small rank:\\

\begin{theorem}Let $\E$ a vector bundle on $\Q_n$ ($n> 3$). If $n$ is odd, $\sS$ the spinor bundle and rank $\E<n-1$, then $\E$
splits if and only if
\begin{enumerate}
\item $H^i_*(\Q_n, \E)=0$ for $2\leq i\leq n-2$ \item $H^1_*(\Q_n,
\E\otimes \sS)=0$.
\end{enumerate}
If $n$ is even, $\sS'$ and $\sS''$ are the two spinor bundles and rank $\E<n$,
then $\E$
splits if and only if
\begin{enumerate}
\item $H^i_*(\Q_n, \E)=0$ for $2\leq i\leq n-2$ \item $H^1_*(\Q_n,
\E\otimes \sS')=H^1_*(\E\otimes\sS'')=0$.
\end{enumerate}
\begin{proof}Let us assume that $\E$ does not split and let us consider a minimal monad for $\E$,

$$ 0 \to \A \xrightarrow{\alpha} \B
\xrightarrow{\beta} \C \to 0.$$ Since $H^i_*(\Q_n, \E)=0$ for $2\leq i\leq n-2$, by (\ref{r1}), $\B$ is an ACM bundle on $Q_n$
and, it has to be isomorphic to a direct sum of line bundles and
spinor bundles twisted by some $\sO(t)$.\\
If $\sS$ is a spinor bundle and $H^1_*(\E\otimes\sS)=0$, from the two sequences   $$ 0 \to \G\otimes\sS \rightarrow \B\otimes\sS
\rightarrow \C\otimes\sS \to 0$$ and $$ 0 \to \A\otimes\sS
\rightarrow \G\otimes\sS \rightarrow \E\otimes\sS \to 0,$$  we can
see that  also $H^1_*(\B\otimes\sS)=0$.\\
Now, in the odd case, since $H^1_*(\sS\otimes\sS)\not=0$ see (\ref{r2}), we can say that no spinor bundle can appear in $\B$. So $\B$ has to
split and this is a contradiction.\\
In the even case, since, according with (\ref{r2}), when $n\equiv 2 \pmod 4$,      $$H^1_*(\sS'\otimes\sS')\not=0$$ and $$H^1_*(\sS''\otimes\sS'')\not=0,$$ or, when $n\equiv 0 \pmod 4$, $$H^1_*(\sS'\otimes\sS'')\not=0,$$ we can say  that
no spinor bundles can appear in $\B$. So $\B$ has to split and this
is a contradiction.

\end{proof}
\end{theorem}

This theorem is the equivalent in $\Q_n$ of the result by Kumar, Peterson and Rao.

\begin{remark} The techniques of this proof are similar to those used by Arrondo and Gra$\tilde{n}$a on the Grassmannian $G(1,4)$ (see \cite{AG}).
\end{remark}

\section{Rank $2$ Bundles without Inner\\ Cohomology}
Let us study more carefully the rank $2$ bundles in $\Q_n$ ($n>3$).\\
In $\Q_4$ by (\ref{r2})  we have that $$H^1_*(\sS'\otimes\sS')=H^1_*(\sS''\otimes\sS'')=0$$ and $$H^1_*(\sS'\otimes\sS'')=\mathbb
C.$$ So from the sequence (see \cite{Ot2}) $$ 0 \to \sS'
\rightarrow \sO_{Q_4}^{\oplus 4} \rightarrow \sS''(1) \to 0,$$ and his dual we see that $$H^2_*(\sS'\otimes\sS')=H^2_*(\sS''\otimes\sS'')=\mathbb
C,$$

 It is then possible  to prove the following theorem:
\begin{theorem}For an  indecomposable rank $2$ bundle $\E$ on $\Q_4$ with\\ $H_*^1(\E)\not=0$ and $H_*^2(\E)=0$, the only possible minimal monad with $\A$ or $\C$ different from zero  is (up to a twist)
\begin{equation} 0 \to \sO \rightarrow
\sS'(1)\oplus\sS''(1) \rightarrow \sO(1) \to 0,
\end{equation}
and such a monad exists.
\begin{proof}First of all in a  minimal monad for $\E$,
$$ 0 \to \A \xrightarrow{\alpha} \B
\xrightarrow{\beta} \C \to 0,$$ $\B$  is an ACM bundle on $\Q_4$;
then it has to be isomorphic to a direct sum of line bundles and
spinor bundles twisted by some $\sO(t)$.\\
 Since $\B$ cannot split at least a spinor bundle must appear.\\
Assume that just one copy of $\sS'$ or one copy of $\sS''$ it appears in $\B$. Since $$\textrm{$\rk$ $\sS''$ $=$ $\rk$ $\sS'$ $=2$}$$ and then $\wedge^2 \sS'$ and $\wedge^2\sS''$ are line bundles, also the
bundle $\wedge^2\B$ is ACM and  the condition
$$H^1_*(\wedge^2\B)\not=0,$$ in (\ref{t2}), is not satisfied.\\
 Assume that  more than one copy of $\sS'$ or more than one copy of $\sS''$  appears in $\B$. Then in the
bundle $\wedge^2\B$, $(\sS'\otimes\sS')(t)$ or
$(\sS''\otimes\sS'')(t)$ appears  and, since
$$H^2_*(\sS'\otimes\sS')=H^2_*(\sS''\otimes\sS'')=\mathbb
C,$$  the condition
$$H^2_*(\wedge^2\B)=0$$ in (\ref{t2}), fails to be satisfied. So $\B$ must contain both $\sS'$ and $\sS''$ 
with some twist and only one copy of each. We can conclude that $\B$  has to be of the form $$ 
 (\bigoplus_i\sO(a_i))\oplus ( \sS'(b))\oplus (\sS''(c)).$$ Let us notice
furthermore that if $H_*^1(\E)$ has more than $1$ generator, $\rk$ $
\C>1$ and $H_*^0( S_2\C)$ has at least $3$ generators.\\
But
$$H_*^1(\wedge^2\B)\backsimeq H_*^1(\sS'\otimes\sS'')=\mathbb C$$ has just $1$ generator and
this is a
contradiction because by (\ref{t2}) $$ \beta_0(H^1_*(\wedge^2\B))\geq
\beta_0(H^0_*(S_2\C)).$$ This means that $\rk$ $\A$ $=$ $\rk$ $\C$ $=1$.\\
At this point the only possible minimal monads are like
$$ 0 \to \sO(-a+c_1(\E)) \rightarrow
\sS'(b)\oplus\sS''(c) \rightarrow \sO(a) \to 0.$$ where $a, b$ and $c$
are integer numbers.\\
Since $\B$ must be isomorphic to $\B^{\vee}(c_1(\E))$ and
$\sS'^{\vee}\backsimeq \sS'(1)$ and $\sS''^{\vee}\backsimeq \sS''(1)$, we
have that $$b=c=\frac{1+c_1(\E)}{2};$$ this means that $c_1(\E)$
must be odd so we can assume $c_1(\E)=-1$ and $b=c=0.$
Now our monad, twisted by $\sO(a+1)$ looks like
$$ 0 \to \sO \xrightarrow{\alpha}
\sS'(a+1)\oplus\sS''(a+1) \rightarrow \sO(2a+1) \to 0$$ and we can
assume $a\geq 0$ because both $\sS'(l)$ and $\sS''(l)$ do have sections
only if $l\geq 1$.\\
 It is possible to
have an injective map $\alpha$ at level of bundles only if
$$c_4(\sS'(a+1)\oplus\sS''(a+1))=c_4(\sS'^{\vee}(a)\oplus{\sS''}^{\vee}(a))=0.$$
 Our goal now is to find the values of $a$ such that this condition is
satisfied.\\
We know (see \cite{Fr}) the intersection ring of $\Q_4$:
$$A^*(\Q_4)=\mathbb Z e_1 \oplus (\mathbb Z e_2\oplus\mathbb Z
e'_2)\oplus\mathbb Z e_3\oplus \mathbb Z e_4.$$ 
We also know that $c_1(\sS'^{\vee})=c_1({\sS''}^{\vee})=1,$
$c_2(\sS'^{\vee})=(1, 0)=e_2$\\
 and $c_2({\sS''}^{\vee})=(0, 1)=e'_2.$
Then

$$c_2 (\sS'^{\vee}(a))=e_2+ a e_1 * (1)e_1 + a e_1 * a e_1 = (1+a+a^2)
e_2 + (a+a^2) e'_2$$ and $$c_2 ({\sS''}^{\vee}(a))=e'_2+ a e_1 *
(1)e_1 + a e_1 * a e_1 = (a+a^2) e_2 + (1+a+a^2) e'_2;$$ so
$$c_4(\sS'^{\vee}(a)\oplus{\sS''}^{\vee}(a))=c_2(\sS'^{\vee}(a)) * c_2({\sS''}^{\vee}(a))=$$
$$=(1+a+a^2)(a+a^2) e_4 + (a+a^2)(1+a+a^2) e_4=2(1+a+a^2)(a+a^2)
e_4$$ This is zero if and only if $a=0$ or $a=-1$ and we can not accept the last case.\\
For $a=0$ we have the claimed monad $$ 0 \to \sO
\xrightarrow{\alpha} \sS'(1)\oplus\sS''(1) \xrightarrow{\beta} \sO(1)
\to 0.$$

We finally want to prove that such a monad exists.\\
We denote by $\sZ_4(1)$ the homology of our monad. We compute $c_1(\sZ_4)=-1$, $c_2(\sZ_4)=(1,1)$ and $H^0(\sZ_4)=0$ and by (\cite{AS} Proposition p. 205) we can conclude that the bundle $\sZ_4$   lies in a sequence
$$ 0 \to \sO \rightarrow \sZ_4(1)
\rightarrow \sI_Y(1) \to 0$$
where $Y$ is the disjoint union of a plane in $\Lambda$ and a plane in $\Lambda'$, the two families of planes in $\Q_4$.\\
We can hence conclude that our monad exists because it is the homology of a well known bundle.
\end{proof}
\end{theorem}
\begin{remark} We can say then that there exist only three  rank $2$ bundles without inner cohomology in $\Q_4$. They are $\sS$, $\sS'$ and  $\sZ_4$ that  is associated, by the Serre correspondence, to two disjoint planes, one in $\Lambda$ and one in $\Lambda'$.
\end{remark}

\begin{corollary} In higher dimension we have:

\begin{enumerate}
\item For an indecomposable rank $2$ bundle $\E$ on $\Q_5$ with
$H_*^2(\E)=0$ and $H_*^3(\E)=0$, the only possible minimal monad with $\A$ or $\C$ not zero  is (up to a twist)
$$ 0 \to \sO \rightarrow
\sS_5(1) \rightarrow \sO(1) \to 0.$$
and such a monad exists.
 \item For $n>5$, no indecomposable bundle of rank $2$ in $\Q_n$ exists with\\
$H^2_*(\E)= ...=H^{n-2}_*(\E)=0$.
\end{enumerate}

\begin{proof}
First of all let us notice that for $n>4$ there is no indecomposable ACM rank $2$ bundle
since the spinor bundles have rank greater than $2$.\\
 Let us then assume
 that $H^1_*(\E)\not=0$ and let us see how many minimal monads it is possible
to find:
\begin{enumerate}
\item
In a  minimal monad for $\E$ in $\Q_5$,
$$ 0 \to \A \xrightarrow{\alpha} \B
\xrightarrow{\beta} \C \to 0,$$ $\B$  is an ACM bundle on $\Q_5$;
then it has to be isomorphic to a direct sum of line bundles and
spinor bundles twisted by some $\sO (t)$,

Moreover, since $H_*^2(\E)=0$ and $H_*^3(\E)=0$,
$\E_{|\Q_4}=\F$ is a bundle with $H^2_*(\F)=0$ and for (\ref{l1}) his 
minimal monad is just the restriction of the  minimal monad for $\E$ $$ 0 \to \A \xrightarrow{\alpha} \B
\xrightarrow{\beta} \C \to 0.$$
For the theorem above, hence, this minimal monad must be
$$ 0 \to \sO \rightarrow
\sS'(1)\oplus\sS''(1) \rightarrow \sO(1) \to 0.$$

Now, since $$\sS_{5_{|\Q_4}}\backsimeq \sS'\oplus\sS'',$$ the only bundle of the form 
$$(\bigoplus_i\sO(a_i))\oplus ( \bigoplus_j\sS_5(b_j))$$ having $\sS'(1)\oplus\sS''(1)$ as restriction on $\Q_4$ is $\sS_5(1)$ and then the claimed
monad $$ 0 \to \sO \xrightarrow{\alpha} \sS_5(1)
\xrightarrow{\beta} \sO(1) \to 0$$ is the only possible.

We finally want to prove that such a monad exists.\\
We denote by $\sZ_5(1)$ the homology of our monad.\\
 We compute $c_1(\sZ_5)=-1$,\\ $c_2(\sZ_5)=1$ and $H^0(\sZ_5)=0$ and by (\cite{Ot4} Main Theorem p. 88) we can conclude that the bundle $\sZ_5$  is a Cayley bundle (see \cite{Ot4} for generalities on Cayley bundles).\\
The bundle $\sZ_5$ appear also in \cite{Ta} and \cite{KPR2}.\\
We can hence conclude that our monad exists because it is the homology of a well known bundle.
 \item
In $\Q_6$ we use the same argument but, since
$\sS'_{6_{|\Q_5}}\backsimeq \sS'_5$
and also
$\sS''_{6_{|\Q_5}}\backsimeq \sS''_5,$
we have two possible minimal monads:
$$ 0 \to \sO \rightarrow
\sS'_6(1) \rightarrow \sO(1) \to 0$$ and $$ 0 \to \sO \rightarrow
\sS''_6(1) \rightarrow \sO(1) \to 0.$$ In both sequences  the condition $$\B\simeq\B^{\vee}(c_1)$$ is not satisfied, since
${\sS'_6}^{\vee}\simeq \sS''_6(1)$ and ${\sS''}_6^{\vee}\simeq \sS'_6(1).$\\ 
So they cannot be the minimal monads of a rank $2$ bundles.\\
We can conclude that  no indecomposable bundle of rank $2$ in
$\Q_6$ exists with $H^2_*(\E)=...=H^{4}_*(\E)=0$ and clearly also in
higher dimension it is not possible to find any bundle without inner cohomology.
\end{enumerate}
\end{proof}
\end{corollary}

As a  conclusion, the Kumar-Peterson-Rao theorem tells us that in $\mathbb P^n$ with $n>3$ there are no rank $2$ bundles without inner cohomology while in $\Q_n$ with $n>3$ there are $4$ of them: precisely $3$ in $\Q_4$ and $1$ in $\Q_5$.\\

It is surprising that this classification of rank $2$ bundle on $\mathbb P^n$ and $\Q_n$ ($n>3$) exactly agrees with the classification by Ancona, Peternell and Wisniewski of rank $2$ Fano bundles (see \cite{APW}).\\
\begin{theorem}[Ancona, Peternell and Wisniewski]
Let $\E$ be a rank $2$ Fano bundle on $\mathbb P^n$ ($n>3$). Then $\E$ splits. \\
Let $\E$ be a rank $2$ Fano bundle on $\Q_n$ ($n>3$). Then either $\E$ splits or:
\begin{enumerate}\item $n=4$ and $\E$ is (up to twist) a spinor bundle or the bundle $\Z_4$.
\item $n=5$ and $\E$ is (up to twist) a Cayley bundle.
\end{enumerate}
\end{theorem}
\begin{corollary} If $\E$ is a rank $2$ bundle on $\mathbb P^n$ and $\Q_n$ ($n>3$), then $$\textrm{$\E$ is a Fano bundle $\Leftrightarrow \E$ is without inner cohomology.}$$
  \end{corollary}

\bibliographystyle{amsplain}

\end{document}